\documentclass[11pt]{article}
\usepackage{amsmath}
\usepackage{amssymb}
\usepackage{amsfonts}
\usepackage{color}
\allowdisplaybreaks

\newtheorem{Theorem}{\textbf{Theorem}}[section]
\newtheorem{Lemma}{\textbf{Lemma}}[section]
\newtheorem{Proposition}{\textbf{Proposition}}[section]
\newtheorem{Corollary}{\textbf{Corollary}}[section]
\newtheorem{Remark}{\textbf{Remark}}[section]
\newtheorem{Example}{\textbf{Example}}[section]
\newtheorem{Definition}{\textbf{Definition}}[section]

\newenvironment{theorem}{\begin{Theorem}$\!\!\!$}{\end{Theorem}}
\newenvironment{lemma}{\begin{Lemma}$\!\!\!$}{\end{Lemma}}

\newenvironment{remark}{\begin{Remark}$\!\!\!$}{\end{Remark}}

\numberwithin{equation}{section}


\begin{document}
\title{Optimal leading term of solutions to wave equations \\
with strong damping terms}
\author{Hironori Michihisa${}^\ast$\\ 
{\small Department of Mathematics, Graduate School of Science, Hiroshima University} \\
{\small Higashi-Hiroshima 739-8526, Japan}
}

\date{%\today
}

\maketitle

\begin{abstract}
We analyze the asymptotic behavior of solutions to wave equations with strong damping terms in $\textbf{R}^n$ $(n\ge1)$,
\[
u_{tt}-\Delta u-\Delta u_t=0, 
\qquad
u(0,x)=u_0(x), 
\quad
u_t(0,x)=u_1(x).
\]
If the initial data belong to suitable weighted $L^1$ spaces, 
lower bounds for the difference between the solutions and the leading terms in the Fourier space are obtained, 
which implies the optimality of expanding methods and some estimates proposed in \cite{Mi} and in this paper. 
\end{abstract}

\footnote[0]{\hspace{-2em} ${}^\ast$Corresponding author.}
\footnote[0]{\hspace{-2em} \textit{Email:} hi.michihisa@gmail.com}
\footnote[0]{\hspace{-2em} 2010 \textit{Mathematics Subject Classification}. 35B40, 35L25, 35L35}
\footnote[0]{\hspace{-2em} \textit{Keywords and Phrases}: Strongly damped wave equation, Asymptotic expansion, Lower bound, Moment condition}

%%%%%%%%%%%%%%%%%%%%%%%%
\section{Introduction}
%\label{sec:1}
%%%%%%%%%%%%%%%%%%%%%%%%
In this paper we consider the Cauchy problem of the solution to the linear strongly damped wave equation 
\begin{equation}
\label{1.1}
\begin{cases}
u_{tt}-\Delta u-\Delta u_t=0, & t>0,\quad x\in\textbf{R}^n, \\
u(0,x)=u_0(x),
\quad
u_t(0,x)=u_1(x),& x\in\textbf{R}^n, \\
\end{cases}
\end{equation}
where $n\ge1$ and $u_0$, $u_1\in L^2(\textbf{R}^n)\cap L^{1,\gamma}(\textbf{R}^n)$. 
The weighted $L^1$ space $L^{1,\gamma}(\textbf{R}^n)$ $(\gamma\ge0)$ is defined by 
\[
L^{1,\gamma}(\textbf{R}^n)
:=\left\{
f\in L^1(\textbf{R}^n):
\|f\|_{1,\gamma}
:=\int_{\textbf{R}^n} (1+|x|)^\gamma |f(x)|\,dx<\infty
\right\}.
\]

Celebrated mathematical results for equation~\eqref{1.1} are $L^p$-$L^q$ decay estimates obtained by Ponce~\cite{P} and Shibata~\cite{S}. 
Later seveval mathematicians have studied wave equations with structural damping terms such as $(-\Delta)^\theta u_t$ and 
exterior domain cases for \eqref{1.1}. 
See e.g., \cite{CDI}, \cite{DE}, \cite{DR-2}, \cite{I}, \cite{IN}, \cite{IT}, \cite{K} and references therein. 
Especially, in \cite{DR}, they have obtained $(L^1\cap L^2)$-$L^2$ and $L^2$-$L^2$ estimates for the solution to \eqref{1.1}. 
They also studied the critical exponent for equation~\eqref{1.1} with nonlinear terms 
but asymptotic profiles were not investigated. 

Now we focus on reviewing known results on the asymptotic behavior of the solution to \eqref{1.1}.  
Recently, Ikehata-Todorova-Yordanov~\cite{ITY} obtained the asymptotic profile in the abstract framework including \eqref{1.1}. 
Limited to equation~\eqref{1.1}, they indicated that the solution of \eqref{1.1} behaves like 
\[
e^{-\frac{t|\xi|^2}{2}}\cos(t|\xi|)\widehat{u_0}
+e^{-\frac{t|\xi|^2}{2}}\frac{\sin(t|\xi|)}{|\xi|}\widehat{u_1},
\qquad 
t\to\infty, 
\]
in the Fourier space. 
However, such diffusion wave properties are not investigated in detail. 
Quite recently, Michihisa~\cite{Mi} have found a way to obtain higher order asymptotic expansions of the solution to \eqref{1.1} in the $L^2$ framework. 
He used the explicit solution formulae in the Fourier space 
and defined some suitable functions generated by evolution operators for \eqref{1.1} to apply the Taylor theorem efficiently. 
It is based on a quite natural and simple idea. 
So the obtained results in \cite{Mi} seems to be optimal, however, the evidence is not shown anywhere. 

When we discuss the optimality of asymptotic expansions, 
appropriate lower bounds need to be shown. 
Concerning this optimality, Ikehata~\cite{I-2} proved that the solution $u=u(t,x)$ of \eqref{1.1} with $n\ge3$ satisfies the following inequalities:
\[
C_1
\left|\int_{\textbf{R}^n} u_1(x)\,dx\right|
t^{-\frac{n}{4}+\frac{1}{2}}
\le \|u(t)\|_2
\le C_2 t^{-\frac{n}{4}+\frac{1}{2}},
\qquad
t\gg1,
\]
where $C_1>0$ is a constant independent of $t$ and the initial data, and $C_2>0$ is a constant independent of $t$. 
Furthermore, the lower dimensional cases $n=1,2$ are also treated by Ikehata-Onodera~\cite{IO}. 
They obtained meaningful inequalities corresponding to the above, 
which says the optimal infinite time blow-up rates are $\sqrt{t}$ $(n=1)$ and $\sqrt{\log t}$ $(n=2)$. 
So in order to obtain $L^2$ bounded solutions 
we have to impose the additional condition that the mass of $u_1$ is zero. 
However, in this case we face the same problem on the decay rate again.
In order to answer this question, in this paper we obtain some asymptotic estimates by 
taking into account the expanding of the initial value, 
which were out of interests in \cite{Mi}, 
and they lead to the optimal lower bound for the $L^2$ difference between the solution to \eqref{1.1} and its leading term.
\\

The rest of this paper is as follows. 
In Section~\ref{sec:2}, we confirm the solution formulae with some notation. 
We also define several functions which are components of asymptotic profiles. 
Main results are stated in Section~\ref{sec:3}. 
Theorem~\ref{thm:3.4} is the most important result in this paper. 
We prepare some lemmas in Section~\ref{sec:4}. 
Expect for Theorem~\ref{thm:3.3}, proofs of theorems are in Section~\ref{sec:5}. 
In Section~\ref{sec:6}, Appendix, 
we give a simpler and complete proof of Theorem~\ref{thm:3.3} which is proved by \cite{IO} for the first time.

%%%%%%%%%%%%%%%%%%%%%%%%
\section{Notation}
\label{sec:2} 
In this paper, $\textbf{N}$ denotes the set of all natural numbers, 
and write $\textbf{N}_0:=\textbf{N}\cup\{0\}$. 
We also write the surface area of the $n$-dimensional unit ball as 
\[
\omega_n
=\int_{|\omega|=1} dS
=\frac{2\pi^\frac{n}{2}}{\Gamma(n/2)}. 
\]
Let us denote the function $\hat{f}$ by the Fourier transform of $f$, 
\begin{equation*}
\hat{f}(\xi)=\int_{\textbf{R}^n}e^{-ix\cdot\xi}f(x)\,dx.
\end{equation*}
%%%
In the Fourier space, the solution $u=u(t,x)$ of \eqref{1.1} is formally expressed by 
\begin{equation}
\label{2.1}
\hat{u}(t,\xi)
=E_0(t,\xi)\widehat{u_0}
+E_1(t,\xi)
\left(
\frac{|\xi|^2}{2}\widehat{u_0}
+\widehat{u_1}
\right), 
\end{equation}
where $E_i$ $(i=0,1)$ are evolution operators given by 
\begin{align*}
E_0(t,\xi)
:=
\begin{cases}
e^{-\frac{t|\xi|^2}{2}}\cos\left(\frac{t|\xi|\sqrt{4-|\xi|^2}}{2}\right), 
& |\xi|\le2, \\
\\
e^{-\frac{t|\xi|^2}{2}}\cosh\left(\frac{t|\xi|\sqrt{|\xi|^2-4}}{2}\right), 
& |\xi|>2,
\end{cases}
\end{align*}
\begin{align*}
E_1(t,\xi)
:=
\begin{cases}
e^{-\frac{t|\xi|^2}{2}}\left[
\sin\left(\frac{t|\xi|\sqrt{4-|\xi|^2}}{2}\right)
\bigg/
\frac{|\xi|\sqrt{4-|\xi|^2}}{2}
\right], 
& |\xi|\le2, \\
\\
e^{-\frac{t|\xi|^2}{2}}\left[
\sinh\left(\frac{t|\xi|\sqrt{|\xi|^2-4}}{2}\right)
\bigg/
\frac{|\xi|\sqrt{|\xi|^2-4}}{2}
\right], 
& |\xi|>2.
\end{cases}
\end{align*}

For example, in \cite{ITY}, it is shown that 
the problem \eqref{1.1} has a unique weak solution 
$u\in C([0,+\infty);H^1(\textbf{R}^n))\cap C^1([0,+\infty);L^2(\textbf{R}^n))$ if $[u_0,u_1]\in H^1(\textbf{R}^n)\times L^2(\textbf{R}^n)$. 
So it is natural to consider $u_0\in H^1(\textbf{R}^n)$ but we treat $u_0\in L^2(\textbf{R}^n)$ in this paper since we can impose additional regularity on the initial datum $u_0$ as far as we need. 

We define the following functions (see \cite{Mi}): 
\begin{align*}
L_0(a,t,\xi) &:=\cos\left(t|\xi|-t|\xi|^2\frac{a}{4+2\sqrt{4-a^2}}\right), \\[5pt]
L_1(a,t,\xi) &:=\sin\left(t|\xi|-t|\xi|^2\frac{a}{4+2\sqrt{4-a^2}}\right)
\bigg/
\frac{|\xi|\sqrt{4-a^2}}{2}.
\end{align*}
Note that 
\[
E_i(t,\xi)=e^{-\frac{t|\xi|^2}{2}}L_i(|\xi|,t,\xi), 
\qquad 
i=0,1. 
\]
Let $k\in\textbf{N}_0$.
Put 
\[
e_i^k(t,\xi)
:=e^{-\frac{t|\xi|^2}{2}}
\frac{1}{k!}\frac{\partial^k L_i}{\partial a^k}(0,t,\xi)
\cdot
|\xi|^k,
\qquad
i=0,1,
\]
and 
\[
m[f]^k(\xi)
:=\sum_{|\alpha|=k} 
\frac{(-1)^{|\alpha|}}{\alpha!}
\left(
\int_{\textbf{R}^n} x^\alpha f(x)\,dx
\right)
(i\xi)^\alpha,
\qquad 
f\in L^{1,k}(\textbf{R}^n). 
\]

%%%%%%%%%%%%%%%%%%%%%%%
\section{Main results}
\label{sec:3}
%%%%%%%%%%%%%%%%%%%%%%%
In \cite{Mi}, he proposed a method for expanding evolution operators $E_i$ $(i=0,1)$ but did not mention any expanding techniques taking into account of the initial data. 
In order to obtain the precise estimate such as \eqref{3.7}, 
we need to prepare more detailed estimates, e.g., \eqref{3.3} and \eqref{3.6}.
\begin{theorem}
\label{thm:3.1}
Let $u_1\in L^{1,\gamma}(\textbf{R}^n)$ with 
\begin{equation}
\label{3.1}
\gamma
\begin{cases}
>\displaystyle{\frac{1}{2}}, & n=1, \\[6pt]
>0, & n=2, \\
\ge0, & n\ge3. 
\end{cases}
\end{equation} 
Then it holds that 
\begin{align}
\label{3.2}
\left\|
E_1(t)\widehat{u_1}
-\sum_{k=0}^{[\gamma]} 
\left(
e_1^k(t)
\sum_{j=0}^{[\gamma]-k} 
m[u_1]^j
\right)
\right\|_{L^2(|\xi|\le1)}
\le C\|u_1\|_{1,\gamma}(1+t)^{-\frac{n}{4}-\frac{\gamma}{2}+\frac{1}{2}},
\qquad
t>0.
\end{align}
Here $C>0$ is a constant independent of $t$ and $u_1$. 
Moreover, it holds that 
\begin{align}
\label{3.3}
\left\|
E_1(t)\widehat{u_1}
-\sum_{k=0}^{[\gamma]} 
\left(
e_1^k(t)
\sum_{j=0}^{[\gamma]-k} 
m[u_1]^j
\right)
\right\|_{L^2(|\xi|\le1)}
=o(t^{-\frac{n}{4}-\frac{\gamma}{2}+\frac{1}{2}}),
\qquad 
t\to\infty.
\end{align}
\end{theorem}

\begin{remark}
\label{rem:3.1}
The condition~\eqref{3.1}, 
which implies $n+2\gamma-2>0$, 
is used to recover the integrability of 
$\sin(t|\xi|)/|\xi|$ near the origin when $n=1,2$. 
Indeed, we can assure that 
\begin{align}
\label{3.4}
\int_{\textbf{R}^n}
|\eta|^{2\gamma-2} e^{-\frac{|\eta|^2}{2}}
\,d\eta 
=\omega_n \int_0^\infty x^{n-3+2\gamma}e^{-\frac{x^2}{2}}\,dx
<\infty.
\end{align}
\end{remark}

\begin{theorem}
\label{thm:3.2}
Let $n\ge1$ and $u_0\in L^{1,\gamma}(\textbf{R}^n)$ with $\gamma\ge0$. 
Then it holds that 
\begin{align}
\label{3.5}
\left\|
E_0(t)\widehat{u_0}
-\sum_{k=0}^{[\gamma]} 
\left(
e_0^k(t)
\sum_{j=0}^{[\gamma]-k} 
m[u_0]^j
\right)
\right\|_{L^2(|\xi|\le1)}
\le C\|u_0\|_{1,\gamma}(1+t)^{-\frac{n}{4}-\frac{\gamma}{2}},
\qquad 
t>0.
\end{align} 
Here $C>0$ is a constant independent of $t$ and $u_0$. 
Moreover, it holds that 
\begin{align}
\label{3.6}
\left\|
E_0(t)\widehat{u_0}
-\sum_{k=0}^{[\gamma]} 
\left(
e_0^k(t)
\sum_{j=0}^{[\gamma]-k} 
m[u_0]^j
\right)
\right\|_{L^2(|\xi|\le1)}
=o(t^{-\frac{n}{4}-\frac{\gamma}{2}}),
\qquad
t\to\infty.
\end{align}
\end{theorem}

The following theorem has been already proved in \cite{I-2} and \cite{IO}. 
We give its simpler proof in Appendix for confirming. 
\begin{theorem}
{\rm \cite{I-2},\cite{IO}}
\label{thm:3.3}
Let $n\ge1$ 
and $u_0$, $u_1\in L^1(\textbf{R}^n)\cap L^2(\textbf{R}^n)$. 
Put 
\[
P_1:=\int_{\textbf{R}^n}u_1(x)\,dx.
\]
Then it holds that
\begin{align*}
C_1^1|P_1|\sqrt{t}
& \le \big\|
\hat{u}(t)
\big\|_2
\le C_2^1\sqrt{t},
& n=1, \\[5pt]
C_1^2|P_1|\sqrt{\log t}
& \le \big\|
\hat{u}(t)
\big\|_2
\le C_2^2\sqrt{\log t},
& n=2, \\[5pt]
C_1^n|P_1|
t^{-\frac{n}{4}+\frac{1}{2}}
& \le \big\|
\hat{u}(t)
\big\|_2
\le C_2^n t^{-\frac{n}{4}+\frac{1}{2}},
& n\ge3, 
\end{align*}
for sufficiently large $t$. 
Here $C_1^n>0$ $(n\ge1)$ are constants independent of $t$ and the initial data, 
and $C_2^n>0$ $(n\ge1)$ are constants independent of $t$. 
\end{theorem}

The goal of this paper is to obtain the following theorem. 
When the mass of $u_1$ is zero, we cannot obtain information on lower bounds from Theorem~\ref{thm:3.3}. 
However, as we can see Theorem~\ref{thm:3.4} below, 
even if 
\[
P_1=\int_{\textbf{R}^n}u_1(x)\,dx=0,
\]
the right-hand side of \eqref{3.7} can be positive as long as one of the quantities 
\[
\int_{\textbf{R}^n}x_j u_1(x)\,dx
\quad 
(j=1,\dots,n)
\qquad
\mbox{and}
\qquad
\int_{\textbf{R}^n} u_0(x)\,dx
\]
is not equal to zero. 
This is one of our advantage of higher order expansions. 
Theorems~\ref{thm:3.3} and \ref{thm:3.4} indicate that moments of initial data are important quantities for precise asymptotic profiles. 
\begin{theorem}
\label{thm:3.4}
Let $n\ge1$ 
and $\hat{u}$ be the function defined by \eqref{2.1} with 
$u_1\in L^{1,1}(\textbf{R}^n)\cap L^2(\textbf{R}^n)$, $u_0\in L^1(\textbf{R}^n)\cap L^2(\textbf{R}^n)$. 
Then it holds that 
\begin{equation}
\label{3.7}
\begin{split}
& \biggr\|
\hat{u}(t)
-\left(
\int_{\textbf{R}^n}u_1(x)\,dx
\right)
e^{-\frac{t|\xi|^2}{2}}\frac{\sin(t|\xi|)}{|\xi|}
\biggr\|_2 \\
& \ge C\sqrt{
\left(
\int_{\textbf{R}^n} u_1(x)\,dx
\right)^2
+\sum_{j=1}^n
\left(
\int_{\textbf{R}^n}x_j u_1(x)\,dx
\right)^2
+\left(
\int_{\textbf{R}^n} u_0(x)\,dx
\right)^2
}
t^{-\frac{n}{4}}
\end{split}
\end{equation}
for sufficiently large $t$. 
Here $C>0$ is a constant independent of $t$. 
\end{theorem}

\begin{remark}
%\label{rem:3.2}
The constant $C>0$ in \eqref{3.7} can be taken not to depend on the initial data if 
\[
\int_{\textbf{R}^n} u_1(x)\,dx=0
\qquad
\mbox{or}
\qquad
\int_{\textbf{R}^n} u_0(x)\,dx=0.
\]
In general, it depends on the ratio 
\[
\left|
\int_{\textbf{R}^n} u_0(x)\,dx
\right|
\biggr/
\left|
\int_{\textbf{R}^n} u_1(x)\,dx
\right|.
\]
See the latter part of the proof of Theorem~\ref{thm:3.4} for details. 
\end{remark}

When we consider the case $u_0\in L^{1,1}(\textbf{R}^n)\cap L^2(\textbf{R}^n)$, $u_1\in L^1(\textbf{R}^n)\cap L^2(\textbf{R}^n)$, 
it follows from \eqref{3.2} with $\gamma=1$ that 
\begin{align*}
\left\|
E_1(t)\widehat{u_1}
-\sum_{k=0}^1 
\left(
e_1^k(t)
\sum_{j=0}^{1-k} 
m[u_1]^j
\right)
\right\|_{L^2(|\xi|\le1)}
\le C\|u_1\|_{1,1}(1+t)^{-\frac{n}{4}},
\qquad
t>0.
\end{align*}
Direct calculation shows 
\begin{align*}
\left\|
\sum_{k=0}^1 
e_1^k(t)
m[u_1]^{1-k}
\right\|_{L^2(|\xi|\le1)}
\le C\|u_1\|_{1,1}(1+t)^{-\frac{n}{4}},
\qquad
t>0.
\end{align*}
So we obtain 
\begin{equation}
\label{3.8}
\biggr\|
\hat{u}(t)
-\left(
\int_{\textbf{R}^n}u_1(x)\,dx
\right)
e^{-\frac{t|\xi|^2}{2}}\frac{\sin(t|\xi|)}{|\xi|}
\biggr\|_2
\le CI(u_1,u_0)(1+t)^{-\frac{n}{4}},
\qquad 
t>0,
\end{equation}
where $I(u_1,u_0):=\|u_1\|_{1,1}+\|u_0\|_1+\|u_1\|_2+\|u_0\|_2$. \\

Inequalities~\eqref{3.7} and \eqref{3.8} imply the optimality of the leading term 
i.e., the zero-th order asymptotic expansion and the decay estimate for the difference between the solution to \eqref{1.1} and the leading term. 
%%%%%%%%%%%%%
\section{Preliminaries}
\label{sec:4}
%%%%%%%%%%%%%
Although, proofs of results in \cite{Mi} were mainly focused on the analysis of $E_0$, 
we deal with $E_1$ intensively. 
The following lemma is a fundamental result for expanding evolution operators. 
\begin{lemma}
{\rm \cite{Mi}}
\label{lem:4.1}
Let $n\ge1$ and $p\in\textbf{N}_0$. 
Then there exists a constant $C>0$ such that 
\begin{align}
\label{4.1}
\left\|
E_1(t)-\sum_{k=0}^p e_1^k(t)
\right\|_{L^2(|\xi|\le1)}
\le C(1+t)^{-\frac{n}{4}-\frac{p}{2}},
\end{align}
\begin{align}
\label{4.2}
\left\|
E_0(t)-\sum_{k=0}^p e_0^k(t)
\right\|_{L^2(|\xi|\le1)}
\le C(1+t)^{-\frac{n}{4}-\frac{p+1}{2}},
\end{align}
for $t>0$. 
\end{lemma}
\textbf{Proof.} 
We reconfirm the proof of \eqref{4.1}. 
Here we consider $\xi\in\textbf{R}^n$ with $|\xi|\le1$. 
By the Taylor theorem we see that 
\begin{align*}
E_1(t,\xi)-\sum_{k=0}^p e_1^k(t,\xi)
& =e^{-\frac{t|\xi|^2}{2}}
\left[
L_1(|\xi|,t,\xi)
-\sum_{k=0}^p\frac{1}{k!}\frac{\partial^k L_1}{\partial a^k}(0,t,\xi)
\cdot
|\xi|^k
\right] \\
& =e^{-\frac{t|\xi|^2}{2}}
\frac{1}{(p+1)!}\frac{\partial^{p+1} L_1}{\partial a^{p+1}}(\tau|\xi|,t,\xi)
\cdot
|\xi|^{p+1}
\end{align*}
for some $0\le\tau\le1$. 
The function 
\[
\frac{a}{4+2\sqrt{4-a^2}}
\]
and its derivatives are all bounded for $0\le a\le1$. 
Thus, for $k\in\textbf{N}_0$, there exists a constant $C>0$ such that 
\begin{align}
\label{4.3}
\left|
\frac{\partial^k L_1}{\partial a^k}(a,t,\xi)
\right|
\le C\frac{e^{-\frac{t|\xi|^2}{2}}}{|\xi|}
\sum_{\ell=0}^k (t|\xi|^2)^\ell
\end{align}
for $0\le a\le1$, $t>0$ and $\xi\in\textbf{R}^n$. 
So we obtain 
\begin{align*}
\left\|
E_1(t)
-\sum_{k=0}^p e_1^k(t)
\right\|_{L^2(|\xi|\le1)}
\le C\sum_{\ell=0}^{p+1}
\left\|
(t|\xi|^2)^\ell 
|\xi|^p
e^{-\frac{t|\xi|^2}{2}}
\right\|_{L^2(|\xi|\le1)} 
\le C(1+t)^{-\frac{n}{4}-\frac{p}{2}}
\end{align*}
for $t>0$. 
One can similarly prove inequality~\eqref{4.2}.
$\Box$
\\

The following lemma is used to expand the Fourier transform of the initial data. 
Lemma~\ref{lem:4.2} is already proved by \cite{IM} but we give a simpler proof than that of \cite{IM}. 
\begin{lemma}
{\rm \cite{IM}}
\label{lem:4.2}
Let $n\ge1$ and $\gamma\ge0$. 
\begin{itemize}
\item[\rm{(i)}] 
It holds that 
\begin{align}
\label{4.4}
\left|
\hat{f}(\xi)
-\sum_{k=0}^{[\gamma]} m[f]^k(\xi)
\right|
\le C|\xi|^\gamma 
\int_{\textbf{R}^n}|x|^\gamma |f(x)|\,dx,
\qquad 
\xi\in\textbf{R}^n, 
\end{align}
for $f\in L^{1,\gamma}(\textbf{R}^n)$. 
Here $C>0$ is a constant independent of $\xi$ and $f$\rm{;}
\item[\rm{(ii)}] 
For $c>0$ and $f\in L^{1,\gamma}(\textbf{R}^n)$, one has 
\begin{align}
\label{4.5}
\lim_{t\to\infty}
t^{\frac{n}{4}+\frac{\gamma}{2}}
\left\|
e^{-ct|\xi|^2}\hat{f}
-\sum_{k=0}^{[\gamma]} m[f]^k e^{-ct|\xi|^2}
\right\|_2
=0.
\end{align}
\end{itemize}
\end{lemma}
\textbf{Proof.} 
First, we prove (i). 
The Taylor theorem yields 
\begin{align*}
e^{-ix\cdot\xi}
-\sum_{|\alpha|\le[\gamma]}
\frac{(-1)^{|\alpha|}}{\alpha!}
x^\alpha (i\xi)^\alpha
=\frac{1}{[\gamma]!}
\int_0^1 (1-\tau)^{[\gamma]}
\frac{d^{[\tau]+1}}{d\tau^{[\gamma]+1}}e^{-i\tau x\cdot\xi}\,d\tau.
\end{align*}
We calculate 
\begin{align*}
\mathcal{R}_{[\gamma]+1}(x,\xi)
& :=\frac{1}{[\gamma]!}
\int_0^1 (1-\tau)^{[\gamma]}
\frac{d^{[\tau]+1}}{d\tau^{[\gamma]+1}}e^{-i\tau x\cdot\xi}\,d\tau \\
& =\frac{(-i)^{[\gamma]+1}}{[\gamma]!}
(x\cdot\xi)^{[\gamma]+1}
\int_0^1 (1-\tau)^{[\gamma]}e^{-i\tau x\cdot\xi}\,d\tau.
\end{align*}
Since 
\[
\sup_{|x||\xi|\ge1}
\frac{\big|\mathcal{R}_{[\gamma]+1}(x,\xi)\big|}{|x|^\gamma |\xi|^\gamma}
\le \sup_{|x||\xi|\ge1}
\left(
\frac{1}{(|x||\xi|)^\gamma}
+\sum_{|\alpha|\le[\gamma]}
\frac{1}{\alpha!}
\frac{1}{(|x||\xi|)^{\gamma-|\alpha|}}
\right)
<\infty,
\]
there exists a constant $C>0$ such that 
\[
\big|\mathcal{R}_{[\gamma]+1}(x,\xi)\big|
\le C|x|^\gamma |\xi|^\gamma 
\]
for $x\in\textbf{R}^n$ and $\xi\in\textbf{R}^n$.
Hence we obtain 
\begin{align*}
& \left|
\hat{f}(\xi)
-\sum_{k=0}^{[\gamma]} m[f]^k(\xi)
\right|
=\left|
\hat{f}(\xi)
-\sum_{|\alpha|\le[\gamma]} 
\frac{(-1)^{|\alpha|}}{\alpha!}
\left(
\int_{\textbf{R}^n}
x^\alpha f(x)\,dx
\right)
(i\xi)^\alpha
\right| \\
& =\left|
\int_{\textbf{R}^n}
\left(
e^{-ix\cdot\xi}
-\sum_{|\alpha|\le[\gamma]}
\frac{(-1)^{|\alpha|}}{\alpha!}
x^\alpha (i\xi)^\alpha
\right)
f(x)\,dx
\right|
=\left|
\int_{\textbf{R}^n}
\mathcal{R}_{[\gamma]+1}(x,\xi)
f(x)\,dx
\right| \\
& \le C|\xi|^\gamma 
\int_{\textbf{R}^n}|x|^\gamma |f(x)|\,dx 
\end{align*}
for $f\in L^{1,\gamma}(\textbf{R}^n)$. 

Next, we check (ii). 
First, we calculate 
\begin{align*}
& \left\|
e^{-ct|\xi|^2}\hat{f}
-\sum_{k=0}^{[\gamma]} m[f]^k e^{-ct|\xi|^2}
\right\|_2^2 
=\int_{\textbf{R}^n}
e^{-2ct|\xi|^2}
\left|
\int_{\textbf{R}^n}
\mathcal{R}_{[\gamma]+1}(x,\xi)
f(x)\,dx
\right|^2
\,d\xi \\
& =t^{-\frac{n}{2}-\gamma}
\int_{\textbf{R}^n}
e^{-2c|\xi|^2}
\left|
t^\frac{\gamma}{2}\int_{\textbf{R}^n}
\mathcal{R}_{[\gamma]+1}\left(x,\frac{\xi}{\sqrt{t}}\right)
f(x)\,dx
\right|^2
\,d\xi.
\end{align*}
Now take any $\xi\in\textbf{R}^n$. 
Then it follows that 
\begin{align*}
\left|
t^\frac{\gamma}{2}\int_{\textbf{R}^n}
\mathcal{R}_{[\gamma]+1}\left(x,\frac{\xi}{\sqrt{t}}\right)
f(x)\,dx
\right|
& \le Ct^\frac{\gamma}{2} \left|\frac{\xi}{\sqrt{t}}\right|^\gamma 
\int_{\textbf{R}^n}|x|^\gamma |f(x)|\,dx \\
& =\left(
C\int_{\textbf{R}^n}|x|^\gamma |f(x)|\,dx
\right)
|\xi|^\gamma, 
\qquad t>0, 
\end{align*}
and, for fixed $x\in\textbf{R}^n$, one has 
\begin{align*}
\left|
t^\frac{\gamma}{2}
\mathcal{R}_{[\gamma]+1}\left(x,\frac{\xi}{\sqrt{t}}\right)
f(x)
\right|
\le Ct^\frac{\gamma}{2}
|x|^\gamma
\left|\frac{\xi}{\sqrt{t}}\right|^\gamma 
|f(x)|
=\biggr(C|\xi|^\gamma\biggr) 
|x|^\gamma |f(x)|, 
\qquad 
t>0, 
\end{align*}
\begin{align*}
& t^\frac{\gamma}{2}
\mathcal{R}_{[\gamma]+1}\left(x,\frac{\xi}{\sqrt{t}}\right) 
=t^{-\frac{[\gamma]+1-\gamma}{2}}
\frac{(-i)^{[\gamma]+1}}{[\gamma]!}
(x\cdot\xi)^{[\gamma]+1}
\int_0^1 (1-\tau)^{[\gamma]}e^{-i\tau x\cdot\frac{\xi}{\sqrt{t}}}\,d\tau 
\longrightarrow 0, 
\qquad
t\to\infty. 
\end{align*}
So the Lebesgue dominanted convergence theorem gives assertion~(ii). 
$\Box$

\begin{remark}
%\label{rem:4.1}
\begin{itemize}
\item[\rm{(A)}] 
Since 
\[
\sum_{|\alpha|\le[\gamma]}
\frac{(-1)^{|\alpha|}}{\alpha!}
x^\alpha (i\xi)^\alpha
=\sum_{k=0}^{[\gamma]}
\frac{1}{k!}
(-ix\cdot\xi)^k,
\]
we have just applied the Taylor theorem for a single-variable function in the proof of assertion~\rm{(i)}. 
\item[\rm{(B)}] 
On \eqref{4.5}, we can say further 
\begin{align*}
\lim_{t\to\infty}
t^{\frac{n}{4}+\frac{\delta}{2}}
\left\|
e^{-ct|\xi|^2}\hat{f}
-\sum_{k=0}^{[\gamma]} m[f]^k e^{-ct|\xi|^2}
\right\|_2
=0
\end{align*}
for any $\delta<[\gamma]+1$.
\end{itemize}
\end{remark}

The following lemma is a generalized version of the results obtained in \cite{I-2} whose proof is based on the double angle formulae and the Riemann-Lebesgue lemma. 
\begin{lemma}
\label{lem:4.3}
\begin{itemize}
\item[\rm{(i)}]
Let $n\ge1$ and let $\gamma$ satisfy \eqref{3.1}. 
Then it follows that 
\begin{align}
\label{4.6}
\int_{|\xi|\le1} 
|\xi|^{2\gamma} e^{-t|\xi|^2}
\left|
\frac{\sin(t|\xi|)}{|\xi|}
\right|^2\,d\xi  
\ge \frac{\omega_n}{4} 
\left(\int_0^1 x^{2\gamma+n-3} e^{-x^2} \,dx \right)
t^{-\frac{n}{2}-\gamma+1}
\end{align}
for sufficiently large $t${\rm ;} 
\item[\rm{(ii)}] Let $n\ge1$ and $\gamma\ge0$. 
Then it follows that 
\begin{align}
\label{4.7}
\int_{|\xi|\le1} 
|\xi|^{2\gamma} e^{-t|\xi|^2}
\left|
\cos(t|\xi|)
\right|^2\,d\xi
\ge\frac{\omega_n}{4} 
\left(
\int_0^1 
x^{2\gamma+n-1} e^{-x^2}\,dx
\right)
t^{-\frac{n}{2}-\gamma}
\end{align}
for sufficiently large $t$.  
\end{itemize}
\end{lemma}

\begin{remark}
\label{rem:4.2}
\begin{itemize}
\item[\rm{[I]}] 
If $n\ge1$ and $\gamma$ satisfies \eqref{3.1}, 
then it follows from \eqref{3.4} that 
\begin{align*}
\int_{|\xi|\le1} 
|\xi|^{2\gamma} e^{-t|\xi|^2}
\left|
\frac{\sin(t|\xi|)}{|\xi|}
\right|^2\,d\xi  
& \le \int_{|\xi|\le1} 
|\xi|^{2\gamma-2} e^{-t|\xi|^2}
\,d\xi \\
& \le e(1+t)^{-\frac{n}{2}-\gamma+1}
\int_{|\eta|\le\sqrt{t}}
|\eta|^{2\gamma-2} e^{-|\eta|^2}
\,d\eta \\
& =e\omega_n 
(1+t)^{-\frac{n}{2}-\gamma+1}
\int_0^{\sqrt{t}} x^{n+2\gamma-3}e^{-x^2}\,dx \\
& \le C(1+t)^{-\frac{n}{2}-\gamma+1},
\qquad
t>0{\rm ;}
\end{align*}
\item[\rm{[II]}] 
If $n\ge1$ and $\gamma\ge0$, 
one has 
\begin{align*}
\int_{|\xi|\le1} 
|\xi|^{2\gamma} e^{-t|\xi|^2}
\left|
\cos(t|\xi|)
\right|^2
\,d\xi 
& \le \int_{|\xi|\le1} 
|\xi|^{2\gamma} e^{-t|\xi|^2}
\,d\xi \\
& \le e(1+t)^{-\frac{n}{2}-\gamma}
\int_{|\eta|\le\sqrt{t}} 
|\eta|^{2\gamma} e^{-|\eta|^2}
\,d\eta \\
& =e\omega_n
(1+t)^{-\frac{n}{2}-\gamma}
\int_0^{\sqrt{t}} 
x^{2\gamma+n-1} e^{-x^2}\,dx \\
& \le C(1+t)^{-\frac{n}{2}-\gamma},
\qquad
t>0.
\end{align*}
\end{itemize}
\end{remark}

\textbf{Proof of Lemma~\ref{lem:4.3}.} 
Let $n\ge1$ and $\gamma$ satisfy \eqref{3.1}. 
Then we have 
\begin{align*}
& \int_{|\xi|\le 1} 
|\xi|^{2\gamma} e^{-t|\xi|^2}
\left|
\frac{\sin(t|\xi|)}{|\xi|}
\right|^2\,d\xi 
=t^{-\frac{n}{2}-\gamma+1}
\int_{|\eta|\le\sqrt{t}} 
|\eta|^{2\gamma} e^{-|\eta|^2}
\left|
\frac{\sin(\sqrt{t}|\eta|)}{|\eta|}
\right|^2\,d\eta \\
& \ge t^{-\frac{n}{2}-\gamma+1}
\int_{|\eta|\le1} 
|\eta|^{2\gamma} e^{-|\eta|^2}
\left|
\frac{\sin(\sqrt{t}|\eta|)}{|\eta|}
\right|^2\,d\eta \\
& =\omega_n t^{-\frac{n}{2}-\gamma+1}
\int_0^1 x^{2\gamma+n-3} e^{-x^2} 
\sin^2(\sqrt{t}x)\,dx \\
& =\frac{\omega_n}{2} 
\left(\int_0^1 x^{2\gamma+n-3} e^{-x^2} \,dx \right)
t^{-\frac{n}{2}-\gamma+1}
-\frac{\omega_n}{2} t^{-\frac{n}{2}-\gamma+1}
\int_0^1 x^{2\gamma+n-3} e^{-x^2} 
\cos(2\sqrt{t}x)\,dx 
\end{align*}
for sufficiently large $t$. 
The assumption~\eqref{3.1} assures that $x^{2\gamma+n-3} e^{-x^2}\in L^1(0,1)$ (see also \eqref{3.4}) and   
the Riemann-Lebesgue lemma implies 
\[
\lim_{t\to\infty}
\int_0^1 x^{2\gamma+n-3} e^{-x^2} 
\cos(2\sqrt{t}x)\,dx
=0.
\]
So we obtain assertion~(i).  
Assertion~(ii) can be proved similarly. 
$\Box$

%%%%%%%%%%%%%%%%%%%
\section{Proofs of main results}
\label{sec:5}
%%%%%%%%%%%%%%%%%%%
\textbf{Proof of Theorem~\ref{thm:3.1}.} 
If $u_1\in L^{1,\gamma}(\textbf{R}^n)$ with $\gamma$ satisfying \eqref{3.1}, 
we see that 
\begin{align*}
& E_1(t,\xi)\widehat{u_1} \\
& =\left[
\sum_{k=0}^{[\gamma]} e_1^k(t,\xi)
+\left(
E_1(t,\xi)
-\sum_{k=0}^{[\gamma]} e_1^k(t,\xi)
\right)
\right]
\left[
\sum_{k=0}^{[\gamma]} m[u_1]^k(\xi)
+
\left(
\widehat{u_1}(\xi)
-\sum_{k=0}^{[\gamma]} m[u_1]^k(\xi)
\right)
\right] \\
& =\left(
\sum_{k=0}^{[\gamma]} e_1^k(t,\xi)
\right)
\left(
\sum_{k=0}^{[\gamma]} m[u_1]^k(\xi)
\right)
\\ 
& \qquad 
+\left(
\sum_{k=0}^{[\gamma]} e_1^k(t,\xi)
\right)
\left(
\widehat{u_1}(\xi)
-\sum_{k=0}^{[\gamma]} m[u_1]^k(\xi)
\right) 
+\left(
E_1(t,\xi)
-\sum_{k=0}^{[\gamma]} e_1^k(t,\xi)
\right)
\left(
\sum_{k=0}^{[\gamma]} m[u_1]^k(\xi)
\right) \\
& \qquad 
+\left(
E_1(t,\xi)
-\sum_{k=0}^{[\gamma]} e_1^k(t,\xi)
\right)
\left(
\widehat{u_1}(\xi)
-\sum_{k=0}^{[\gamma]} m[u_1]^k(\xi)
\right).
\end{align*}
From \eqref{4.3}, 
there exists a constant $C>0$ such that 
\[
\left|
\sum_{k=0}^{[\gamma]} e_1^k(t,\xi)
\right|
\le C\frac{e^{-\frac{t|\xi|^2}{2}}}{|\xi|}
\sum_{\ell=0}^{[\gamma]} (t|\xi|^2)^\ell
\]
for $t>0$ and $\xi\in\textbf{R}^n$ with $|\xi|\le1$.  
Thus we have 
\begin{align*}
\left\|
\left(
\sum_{k=0}^{[\gamma]} e_1^k(t)
\right)
\left(
\widehat{u_1}
-\sum_{k=0}^{[\gamma]} m[u_1]^k
\right)
\right\|_{L^2(|\xi|\le1)}
& \le C\|u_1\|_{1,\gamma}
\sum_{\ell=0}^{[\gamma]} 
\left\|
(t|\xi|^2)^\ell
|\xi|^{\gamma-1} 
e^{-\frac{t|\xi|^2}{2}}
\right\|_{L^2(|\xi|\le1)} \\
& \le C\|u_1\|_{1,\gamma}
(1+t)^{-\frac{n}{4}-\frac{\gamma}{2}+\frac{1}{2}},
\qquad
t>0,
\end{align*}
with the aid of \eqref{4.4} and [I] in Remark~\ref{rem:4.2}. 
It follows from \eqref{4.1} that 
\[
\left\|
\left(
E_1(t)
-\sum_{k=0}^{[\gamma]} e_1^k(t)
\right)
\left(
\sum_{k=0}^{[\gamma]} m[u_1]^k
\right)
\right\|_{L^2(|\xi|\le1)}
\le C\|u_1\|_{1,[\gamma]}
(1+t)^{-\frac{n}{4}-\frac{[\gamma]}{2}},
\qquad
t>0.
\]
The last term is estimated by 
\begin{align*}
\left\|
\left(
E_1(t)
-\sum_{k=0}^{[\gamma]} e_1^k(t)
\right)
\left(
\widehat{u_1}
-\sum_{k=0}^{[\gamma]} m[u_1]^k
\right)
\right\|_{L^2(|\xi|\le1)}
\le C\|u_1\|_{1,\gamma}
(1+t)^{-\frac{n}{4}-\frac{[\gamma]+\gamma}{2}},
\qquad
t>0.
\end{align*}
Furthermore, one has 
\begin{align*}
& \left(
\sum_{k=0}^{[\gamma]} e_1^k(t,\xi)
\right)
\left(
\sum_{k=0}^{[\gamma]} m[u_1]^k(\xi)
\right) \\
& \qquad
=\sum_{k=0}^{[\gamma]} 
\left(
e_1^k(t,\xi)
\sum_{j=0}^{[\gamma]-k} 
m[u_1]^j(\xi)
\right)
+\sum_{k=1}^{[\gamma]}
\left(
e_1^k(t,\xi)
\sum_{j=[\gamma]-k+1}^{[\gamma]} 
m[u_1]^j(\xi)
\right). 
\end{align*}
If $[\gamma]=0$, then the second term of the right-hand side does not appear. 
So we consider the case $\gamma\ge1$. 
In this case it follows from \eqref{4.3} that 
\begin{align*}
\left\|
\sum_{k=1}^{[\gamma]}
\left(
e_1^k(t)
\sum_{j=[\gamma]-k+1}^{[\gamma]} 
m[u_1]^j
\right)
\right\|_{L^2(|\xi|\le1)} 
& \le C\|u_1\|_{1,[\gamma]}
\sum_{k=1}^{[\gamma]} 
\sum_{\ell=0}^k 
\sum_{j=[\gamma]-k+1}^{[\gamma]} 
\left\|
(t|\xi|^2)^\ell
|\xi|^{k+j-1}
e^{-\frac{t|\xi|^2}{2}}
\right\|_{L^2(|\xi|\le1)} \\
& \le C\|u_1\|_{1,[\gamma]}
(1+t)^{-\frac{n}{4}-\frac{[\gamma]}{2}},
\qquad
t>0.
\end{align*}
Thus we obtain \eqref{3.2}. 

In order to prove \eqref{3.3}, 
it suffices to check 
\begin{align*}
\left\|
\left(
\sum_{k=0}^{[\gamma]} e_1^k(t)
\right)
\left(
\widehat{u_1}
-\sum_{k=0}^{[\gamma]} m[u_1]^k
\right)
\right\|_{L^2(|\xi|\le1)}
=o(t^{-\frac{n}{4}-\frac{\gamma}{2}+\frac{1}{2}}), 
\qquad 
t\to\infty. 
\end{align*}
We calculate 
\begin{align*}
& \left\|
\left(
\sum_{k=0}^{[\gamma]} e_1^k(t)
\right)
\left(
\widehat{u_1}
-\sum_{k=0}^{[\gamma]} m[u_1]^k
\right)
\right\|_{L^2(|\xi|\le1)} 
\le C\sum_{k=0}^{[\gamma]} 
\left\|
(t|\xi|^2)^k
\frac{e^{-\frac{t|\xi|^2}{2}}}{|\xi|}
\left(
\widehat{u_1}
-\sum_{k=0}^{[\gamma]} m[u_1]^k
\right)
\right\|_{L^2(|\xi|\le1)} \\
& \le C\sum_{k=0}^{[\gamma]} 
\sup_{a\ge0}\biggr(a^k e^{-\frac{a}{4}}\biggr)
\left\|
\frac{e^{-\frac{t|\xi|^2}{4}}}{|\xi|}
\left(
\widehat{u_1}
-\sum_{k=0}^{[\gamma]} m[u_1]^k
\right)
\right\|_{L^2(|\xi|\le1)} \\
& \le Ct^{-\frac{n}{4}-\frac{\gamma}{2}+\frac{1}{2}}
\left(
\int_{\textbf{R}^n}
\frac{e^{-\frac{|\eta|^2}{2}}}{|\eta|^2} 
\left|
t^{\frac{\gamma}{2}}
\int_{\textbf{R}^n}
\mathcal{R}_1\left(x,\frac{\eta}{\sqrt{t}}\right) u_1(x)
\,dx
\right|^2\,d\eta
\right)^\frac{1}{2}. 
\end{align*}
In the proof of Lemma~\ref{lem:4.2}, 
we have already shown that 
\begin{align*}
\left|
t^{\frac{\gamma}{2}}
\int_{\textbf{R}^n}
\mathcal{R}_1\left(x,\frac{\eta}{\sqrt{t}}\right) u_1(x)
\,dx
\right|
\le C\|u_1\|_{1,\gamma} |\eta|^\gamma,
\qquad
t>0, 
\end{align*}
for each $\eta\in\textbf{R}^n$.  
Thus the Lebesgue dominanted convergence theorem yields 
\begin{align*}
\lim_{t\to\infty}
\int_{\textbf{R}^n}
\frac{e^{-\frac{|\eta|^2}{2}}}{|\eta|^2} 
\left|
t^{\frac{\gamma}{2}}
\int_{\textbf{R}^n}
\mathcal{R}_1\left(x,\frac{\eta}{\sqrt{t}}\right) u_1(x)
\,dx
\right|^2\,d\eta
=0
\end{align*}
with the aid of \eqref{3.4}. 
So we obtain \eqref{3.3} and the proof is now complete. 
$\Box$
\\

The proof of Theorem~\ref{thm:3.2} is similar to that of Theorem~\ref{thm:3.1}, 
which is simpler. 
Hence we omit the proofs of \eqref{3.5} and \eqref{3.6}. \\

\textbf{Proof of Theorem~\ref{thm:3.4}.} 
It follows from 
\eqref{3.3} with $\gamma=1$ and 
\eqref{3.6} with $\gamma=0$ that 
\begin{align*}
& \biggr\|
\hat{u}(t)
-\left(
\int_{\textbf{R}^n}u_1(x)\,dx
\right)
e^{-\frac{t|\xi|^2}{2}}\frac{\sin(t|\xi|)}{|\xi|}
\biggr\|_2
=\biggr\|
\hat{u}(t)
-e_1^0(t)m[u_1]^0
\biggr\|_2 
\ge \biggr\|
\hat{u}(t)
-e_1^0(t)m[u_1]^0
\biggr\|_{L^2(|\xi|\le1)} 
\\
& \ge \biggr\|
\sum_{k=0}^1
e_1^k(t)
m[u_1]^{1-k} 
+e_0^0(t)m[u_0]^0
\biggr\|_{L^2(|\xi|\le1)} 
-\left\|
E_1(t)\widehat{u_1}
-\sum_{k=0}^1 \left(
e_1^k(t)
\sum_{j=0}^{1-k}
m[u_1]^j
\right)
\right\|_{L^2(|\xi|\le1)} \\
& \qquad\qquad 
-\biggr\|
E_0(t)\widehat{u_0}
-e_0^0(t)m[u_0]^0
\biggr\|_{L^2(|\xi|\le1)}
-\left\|
\frac{|\xi|^2}{2}E_1(t)\widehat{u_1}
\right\|_{L^2(|\xi|\le1)} \\
& \ge \biggr\|
\sum_{k=0}^1
e_1^k(t)
m[u_1]^{1-k}
+e_0^0(t)m[u_0]^0
\biggr\|_{L^2(|\xi|\le1)} 
-o(t^{-\frac{n}{4}})
-o(t^{-\frac{n}{4}})
-O(t^{-\frac{n}{4}-\frac{1}{2}}),
\qquad
t\to\infty.
\end{align*}
Here, the following estimate is just used:
\begin{align}
\label{5.1}
\left\|
\frac{|\xi|^2}{2}E_1(t)
\right\|_{L^2(|\xi|\le1)}^2 
\le C\int_{|\xi|\le1} 
|\xi|^2 e^{-t|\xi|^2}\,d\xi 
\le C(1+t)^{-\frac{n}{2}-1},
\qquad
t>0. 
\end{align}
Recall that 
\begin{align*}
& \sum_{k=0}^1
e_1^k(t,\xi)
m[u_1]^{1-k}(\xi) 
=e_1^0(t,\xi)
m[u_1]^1(\xi)
+e_1^1(t,\xi)
m[u_1]^0(\xi) \\
& =-i \sum_{|\alpha|=1}
\left(
\int_{\textbf{R}^n} x^\alpha u_1(x)\,dx
\right)
e^{-\frac{t|\xi|^2}{2}}\frac{\sin(t|\xi|)}{|\xi|}
\xi^{\alpha} 
-\frac{1}{8}\left(
\int_{\textbf{R}^n} u_1(x)\,dx
\right)
t|\xi|^2 
e^{-\frac{t|\xi|^2}{2}}\cos(t|\xi|),
\end{align*}
\[
e_0^0(t,\xi)m[u_0]^0(\xi)
=\left(
\int_{\textbf{R}^n} u_0(x)\,dx
\right)
e^{-\frac{t|\xi|^2}{2}}\cos(t|\xi|),
\]
and so we have 
\begin{align}
\notag
& \left\|
\sum_{k=0}^1
e_1^k(t)
m[u_1]^{1-k}
+e_0^0(t,\xi)m[u_0]^0 
\right\|_{L^2(|\xi|\le1)}^2 \\
\notag
& =\int_{|\xi|\le1}
e^{-t|\xi|^2}
\left|
i \sum_{|\alpha|=1}
\left(
\int_{\textbf{R}^n} x^\alpha u_1(x)\,dx
\right)
\frac{\sin(t|\xi|)}{|\xi|}
\xi^{\alpha} \right.\\
\notag
& \qquad \left. 
+\frac{1}{8}\left(
\int_{\textbf{R}^n} u_1(x)\,dx
\right)
t|\xi|^2 
\cos(t|\xi|)
-\left(
\int_{\textbf{R}^n} u_0(x)\,dx
\right)
\cos(t|\xi|)
\right|^2
\,d\xi \\
\notag
& =\int_{|\xi|\le1}
e^{-t|\xi|^2}
\left\{
\left[
\sum_{|\alpha|=1}
\left(
\int_{\textbf{R}^n} x^\alpha u_1(x)\,dx
\right)
\frac{\sin(t|\xi|)}{|\xi|}
\xi^{\alpha}
\right]^2 \right. \\
\notag
& \qquad \left.
+\left[
\frac{1}{8}\left(
\int_{\textbf{R}^n} u_1(x)\,dx
\right)
t|\xi|^2 
-\left(
\int_{\textbf{R}^n} u_0(x)\,dx
\right)
\right]^2
|\cos(t|\xi|)|^2
\right\}\,d\xi \\
\notag
& =\sum_{j=1}^n
\left(
\int_{\textbf{R}^n} x_j u_1(x)\,dx
\right)^2
\int_{|\xi|\le1}
e^{-t|\xi|^2}
\left|
\frac{\sin(t|\xi|)}{|\xi|}
\xi_j
\right|^2
\,d\xi \\
\label{5.2}
& \qquad 
+2\sum_{1\le j<k\le n}
\left(
\int_{\textbf{R}^n} x_j u_1(x)\,dx
\right)
\left(
\int_{\textbf{R}^n} x_k u_1(x)\,dx
\right)
\int_{|\xi|\le1}
e^{-t|\xi|^2}
\left|
\frac{\sin(t|\xi|)}{|\xi|}
\right|^2
\xi_j \xi_k
\,d\xi \\
\notag
& \qquad 
+\int_{|\xi|\le1}
e^{-t|\xi|^2}
\left[
\frac{1}{8}\left(
\int_{\textbf{R}^n} u_1(x)\,dx
\right)
t|\xi|^2 
-\left(
\int_{\textbf{R}^n} u_0(x)\,dx
\right)
\right]^2
|\cos(t|\xi|)|^2
\,d\xi.
\end{align}
Now we calculate the first term. 
For all $j=1,\dots,n$, it follows from \eqref{4.6} that 
\begin{align*}
& \int_{|\xi|\le1}
e^{-t|\xi|^2}
\left|
\frac{\sin(t|\xi|)}{|\xi|}
\xi_j
\right|^2
\,d\xi
=\int_{|\xi|\le1}
e^{-t|\xi|^2}
\left|
\frac{\sin(t|\xi|)}{|\xi|}
\right|^2
\xi_1^2
\,d\xi \\
& =\frac{1}{n}\int_{|\xi|\le1}
|\xi|^2
e^{-t|\xi|^2}
\left|
\frac{\sin(t|\xi|)}{|\xi|}
\right|^2
\,d\xi 
\ge \frac{\omega_n}{4n} 
\left(\int_0^1 x^{n-1} e^{-x^2} \,dx \right)
t^{-\frac{n}{2}}
\end{align*}
for sufficiently large $t$.

We can ignore the second term on the right-hand side of \eqref{5.2}. 
This is just because this term never appears if $n=1$, 
and it also holds that 
\[
\int_{|\xi|\le1}
e^{-t|\xi|^2}
\left|
\frac{\sin(t|\xi|)}{|\xi|}
\right|^2
\xi_j \xi_k
\,d\xi
=0,
\qquad
1\le j<k\le n,
\]
when $n\ge2$.
So let us estimate the integral
\begin{align*}
I
& :=\int_{|\xi|\le1}
\left(
\frac{P_1}{8} t|\xi|^2 
-P_0
\right)^2 
e^{-t|\xi|^2}
|\cos(t|\xi|)|^2
\,d\xi \\
& =t^{-\frac{n}{2}}
\int_{|\eta|\le\sqrt{t}}
\left(
\frac{P_1}{8} |\eta|^2 
-P_0
\right)^2 
e^{-|\eta|^2}
|\cos(\sqrt{t}|\eta|)|^2
\,d\eta,
\end{align*}
where 
\[
P_j:=\int_{\textbf{R}^n} u_j(x)\,dx,
\qquad
j=0,1.
\]
If $P_1=0$, then it follows from \eqref{4.7} that 
\begin{align*}
I
& =P_0^2 
\int_{|\xi|\le1}
e^{-t|\xi|^2}
|\cos(t|\xi|)|^2
\,d\xi 
\ge P_0^2
\frac{\omega_n}{4} 
\left(
\int_0^1 
x^{n-1} e^{-x^2}\,dx
\right)
t^{-\frac{n}{2}}, 
\qquad 
t\gg1. 
\end{align*}
Conversely, in the case $P_0=0$, from \eqref{4.7}, one has 
\begin{align*}
I
& =\frac{P_1^2 t^2}{64}
\int_{|\xi|\le1}
|\xi|^4 e^{-t|\xi|^2}
|\cos(t|\xi|)|^2
\,d\xi 
\ge P_1^2
\frac{\omega_n}{256} 
\left(
\int_0^1 
x^{n+3} e^{-x^2}\,dx
\right)
t^{-\frac{n}{2}}, 
\qquad 
t\gg1. 
\end{align*}
Finally, we deal with the case $P_1\not=0$ and $P_0\not=0$. 
Let $\delta>0$ be an arbitrary real number. 
Then, one has 
\begin{align*}
I \ge t^{-\frac{n}{2}}
\int_{\delta\le|\eta|\le2\delta}
\left(
\frac{P_1}{8} |\eta|^2 
-P_0
\right)^2 
e^{-|\eta|^2}
|\cos(\sqrt{t}|\eta|)|^2
\,d\eta 
\end{align*}
for $t\ge4\delta^2>0$. 
Now we choose 
\[
\delta:=4\sqrt{\frac{|P_0|}{|P_1|}}>0.
\]
If $|\eta|\ge\delta$, then 
\[
\left|
\frac{P_1}{8} |\eta|^2 
-P_0
\right|
\ge \frac{|P_1|}{8} |\eta|^2
-|P_0|
\ge \frac{|P_1|}{16}|\eta|^2
\ge |P_0|.
\]
So we can estimate $I$ in two ways: 
\begin{align*}
I & \ge t^{-\frac{n}{2}}
P_0^2
\int_{\delta\le|\eta|\le2\delta} 
e^{-|\eta|^2}
|\cos(\sqrt{t}|\eta|)|^2
\,d\eta \\
& \ge 
\frac{\omega_n}{4}
\left(
\int_{\delta}^{2\delta}
x^{n-1}e^{-x^2}\,dx
\right)
P_0^2 
t^{-\frac{n}{2}},
\end{align*}
\begin{align*}
I & \ge t^{-\frac{n}{2}}
\left(
\frac{P_1}{16}
\right)^2
\int_{\delta\le|\eta|\le2\delta} 
|\eta|^4 e^{-|\eta|^2}
|\cos(\sqrt{t}|\eta|)|^2
\,d\eta \\
& \ge \frac{\omega_n}{1024}
\left(
\int_{\delta}^{2\delta}
x^{n+3}e^{-x^2}\,dx
\right)
P_1^2
t^{-\frac{n}{2}},
\end{align*}
for sufficiently large $t\ge4\delta^2$. 
Therefore we obtain the desired statement of Theorem~\ref{thm:3.4}. 
$\Box$

%%%%%%%%%%%
\section{Appendix}
\label{sec:6}
%%%%%%%%%%%
In this section we give a simplified proof of Theorem~\ref{thm:3.3}. 
Theorem~\ref{thm:3.3} is already proved by \cite{I-2} in the case $n\ge3$. 
Later Theorem~\ref{thm:3.3} with $n=1,2$ is obtained by \cite{IO}. 
In the lower dimensional cases we need some careful calculation on the integral: 
\begin{align*}
\int_{\textbf{R}^n} e^{-t|\xi|^2}
\left|
\frac{\sin(t|\xi|)}{|\xi|}
\right|^2\,d\xi 
=t^{1-\frac{n}{2}}
\int_{\textbf{R}^n} e^{-|\eta|^2}
\left|
\frac{\sin(\sqrt{t}|\eta|)}{|\eta|}
\right|^2\,d\eta. 
\end{align*}
In \cite{I-2}, the following lower bound is obtained 
\begin{align*}
\int_{\textbf{R}^n} e^{-t|\xi|^2}
\left|
\frac{\sin(t|\xi|)}{|\xi|}
\right|^2\,d\xi 
\ge C_n t^{-n+2}, 
\qquad
t\gg1, 
\end{align*}
for all $n\ge1$. 
So the left-hand side of inequality \eqref{6.1} is already proved,  
but \eqref{6.1} was proved in \cite{IO} for the first time in the full version. 
\begin{lemma}
{\rm \cite{I-2}, \cite{IO}}
\label{lem:6.1}
There exist constants $0<C_1^1\le C_2^1<\infty$ such that 
\begin{align}
\label{6.1}
C_1^1 t
\le \int_{-\infty}^\infty e^{-t|\xi|^2}
\left|
\frac{\sin(t|\xi|)}{|\xi|}
\right|^2\,d\xi
\le C_2^1 t
\end{align}
for sufficiently large $t$. 
\end{lemma}
\textbf{Proof.} 
First, we have 
\begin{align*}
& \int_{-\infty}^\infty e^{-t|\xi|^2}
\left|
\frac{\sin(t|\xi|)}{|\xi|}
\right|^2\,d\xi 
=2\sqrt{t}\int_0^\infty 
e^{-\eta^2} 
\left|
\frac{\sin(\sqrt{t}\eta)}{\eta}
\right|^2\,d\eta. 
\end{align*}
Since 
\[
\frac{\theta}{2}\le \sin\theta\le \theta, 
\qquad 
0\le\theta\le1, 
\]
one has 
\[
\frac{\sqrt{t}\eta}{2}\le \sin(\sqrt{t}\eta)\le \sqrt{t}\eta,
\qquad
0\le\eta\le\frac{1}{\sqrt{t}},
\]
for $t\ge1$. 
Thus it follows that 
\begin{align*}
\int_0^{1/\sqrt{t}} 
e^{-\eta^2}
\left|
\frac{\sin(\sqrt{t}\eta)}{\eta}
\right|^2\,d\eta 
\le t\int_0^{1/\sqrt{t}}
e^{-\eta^2}\,d\eta 
\le t\cdot\frac{1}{\sqrt{t}}
=\sqrt{t}
\end{align*}
and 
\begin{align*}
\int_0^{1/\sqrt{t}} 
e^{-\eta^2}
\left|
\frac{\sin(\sqrt{t}\eta)}{\eta}
\right|^2\,d\eta 
& \ge \frac{t}{4}\int_0^{1/\sqrt{t}}
e^{-\eta^2}\,d\eta 
\ge \frac{t}{4}
\int_0^{1/\sqrt{t}} 
e^{-\eta}\,d\eta 
\ge \frac{t}{4}
\int_{1/2\sqrt{t}}^{1/\sqrt{t}} 
e^{-\eta}\,d\eta \\
& \ge \frac{t}{4}
\cdot
\frac{1}{2\sqrt{t}}
e^{-\frac{1}{\sqrt{t}}}
\ge \frac{\sqrt{t}}{8e}
\end{align*}
for $t\ge1$. 
Next, we see that 
\begin{align*}
\int_{1/\sqrt{t}}^\infty 
e^{-\eta^2}
\left|
\frac{\sin(\sqrt{t}\eta)}{\eta}
\right|^2\,d\eta 
& \le \int_{1/\sqrt{t}}^\infty 
\frac{e^{-\eta^2}}{\eta^2}\,d\eta
=\int_{1/\sqrt{t}}^\infty
\left(
-\frac{1}{\eta}
\right)'
e^{-\eta^2}
\,d\eta \\
& =\left[
-\frac{1}{\eta}e^{-\eta^2}
\right]_{1/\sqrt{t}}^\infty
-2\int_{1/\sqrt{t}}^\infty
e^{-\eta^2}
\,d\eta \\
& \le \sqrt{t}e^{-1/t}
\le \sqrt{t}.
\end{align*}
On the other hand, since  
\[
\sin\theta\ge\frac{1}{2}, 
\qquad
1\le\theta\le\frac{\pi}{2},
\]
we have 
\[
\sin(\sqrt{t}\eta)\ge\frac{1}{2},
\qquad
\frac{1}{\sqrt{t}}\le\eta\le\frac{\pi}{2},
\]
when $t\ge1$. 
Thus one has 
\begin{align*}
& \int_{1/\sqrt{t}}^\infty 
e^{-\eta^2}
\left|
\frac{\sin(\sqrt{t}\eta)}{\eta}
\right|^2\,d\eta 
\ge \frac{1}{4}
\int_{1/\sqrt{t}}^{\pi/2}
\frac{e^{-\eta^2}}{\eta^2}\,d\eta 
=\frac{1}{4}\left[
-\frac{1}{\eta}e^{-\eta^2}
\right]_{1/\sqrt{t}}^{\pi/2}
-\frac{1}{2}\int_{1/\sqrt{t}}^{\pi/2}
e^{-\eta^2}
\,d\eta \\
& \ge \frac{1}{4}\biggr(
\sqrt{t}e^{-1/t}
-\frac{2}{\pi}e^{-\pi^2/4}
\biggr)
-\frac{\sqrt{\pi}}{2} 
\ge \frac{1}{4e}\sqrt{t}
-\frac{1}{2\pi}e^{-\pi^2/4}
-\frac{\sqrt{\pi}}{2}.
\end{align*}
Therefore we obtain the lemma. 
$\Box$
\\

As we can easily predict, the two-dimensional case is the most difficult. 
\begin{lemma}
{\rm \cite{IO}}
\label{lem:6.2}
There exist constants $0<C_1^2\le C_2^2<\infty$ such that 
\begin{align}
\label{6.2}
C_1^2 \log t
\le \int_{\textbf{R}^2} e^{-t|\xi|^2}
\left|
\frac{\sin(t|\xi|)}{|\xi|}
\right|^2\,d\xi
\le C_2^2 \log t
\end{align}
for sufficiently large $t$. 
\end{lemma}
\textbf{Proof.} 
Similarly, we have 
\begin{align*}
\int_{|\eta|\le1/\sqrt{t}} 
e^{-|\eta|^2}
\left|
\frac{\sin(\sqrt{t}|\eta|)}{|\eta|}
\right|^2\,d\eta 
\le t\int_{|\eta|\le1/\sqrt{t}} 
e^{-|\eta|^2}\,d\eta 
\le \pi,
\qquad
t>0,
\end{align*}
and since $\sqrt{t}|\eta|\le1$, one has 
\begin{align*}
& \int_{|\eta|\le1/\sqrt{t}} 
e^{-|\eta|^2}
\left|
\frac{\sin(\sqrt{t}|\eta|)}{|\eta|}
\right|^2\,d\eta 
\ge \frac{t}{4}
\int_{|\eta|\le1/\sqrt{t}} 
e^{-|\eta|^2}\,d\eta \\
& \ge \frac{\pi}{2}t 
\int_0^{1/\sqrt{t}} r e^{-r^2}\, dr 
=\frac{\pi}{4}t(1-e^{-1/t}) 
=\frac{\pi}{4}+O\left(\frac{1}{t}\right), 
\qquad 
t\to\infty.
\end{align*}
Here we note that 
\[
1-e^{-1/t}
=\frac{1}{t}
-\sum_{k=2}^\infty 
\frac{1}{k!} \left(-\frac{1}{t}\right)^k
\ge \frac{1}{t}
-\frac{1}{t^2}\sum_{k=2}^\infty 
\frac{1}{k!},
\qquad 
t\ge1.
\]
Next, we see that 
\begin{align*}
& \int_{|\eta|\ge1/\sqrt{t}}
e^{-|\eta|^2}
\left|
\frac{\sin(\sqrt{t}|\eta|)}{|\eta|}
\right|^2\,d\eta 
\le \int_{|\eta|\ge1/\sqrt{t}}
\frac{e^{-|\eta|^2}}{|\eta|^2} \,d\eta \\
& =2\pi \int_{1/\sqrt{t}}^\infty
\frac{e^{-r^2}}{r}\,dr
=2\pi \int_{1/\sqrt{t}}^\infty
(\log r)' e^{-r^2}\,dr \\
& =2\pi\biggr[(\log r)e^{-r^2}\biggr]_{1/\sqrt{t}}^\infty
+4\pi \int_{1/\sqrt{t}}^\infty
(\log r)r e^{-r^2}\,dr \\
& =\pi (\log t)e^{-1/t}
+4\pi \int_{1/\sqrt{t}}^\infty
(\log r)r e^{-r^2}\,dr.
\end{align*}
On the other hand, since  
\[
0\le (\log t)e^{-1/t}\le \log t,
\qquad
t\ge1,
\]
\begin{align*}
\int_0^\infty
(\log r)r e^{-r^2}\,dr
\le \int_0^\infty 
|\log r| r e^{-r^2}\,dr
=:C_0
<\infty, 
\qquad
t>0,
\end{align*}
we arrive at 
\begin{align*}
\int_{|\eta|\ge1/\sqrt{t}}
e^{-|\eta|^2}
\left|
\frac{\sin(\sqrt{t}|\eta|)}{|\eta|}
\right|^2\,d\eta 
\le \int_{|\eta|\ge1/\sqrt{t}}
\frac{e^{-|\eta|^2}}{|\eta|^2} \,d\eta 
\le \pi\log t
+4\pi C_0,
\qquad
t\ge1. 
\end{align*}
Finally, we calculate 
\begin{align*}
& \int_{|\eta|\ge1/\sqrt{t}}
e^{-|\eta|^2}
\left|
\frac{\sin(\sqrt{t}|\eta|)}{|\eta|}
\right|^2\,d\eta 
\ge \frac{1}{4} 
\int_{1/\sqrt{t}\le|\eta|\le\pi/2}
\frac{e^{-|\eta|^2}}{|\eta|^2}\,d\eta \\
& =\frac{\pi}{2}\int_{1/\sqrt{t}}^{\pi/2}
\frac{e^{-r^2}}{r}\,dr 
=\frac{\pi}{2}
\biggr[(\log r)e^{-r^2}\biggr]_{1/\sqrt{t}}^{\pi/2}
+\pi \int_{1/\sqrt{t}}^{\pi/2}
(\log r)r e^{-r^2}\,dr \\
& \ge \frac{\pi}{2}
\left[
\log\left(\frac{\pi}{2}\right)e^{-\pi^2/4}
-\log\left(\frac{1}{\sqrt{t}}\right)e^{-1/t}
\right]
-\pi C_0 \\
& =\frac{\pi}{4e}\log t
+\frac{\pi}{2}\log\left(\frac{\pi}{2}\right)e^{-\pi^2/4}
-\pi C_0,
\qquad
t\ge1.  
\end{align*}
The proof is now complete. 
$\Box$

\begin{lemma}
{\rm \cite{I-2}}
\label{lem:6.3}
Let $n\ge3$. 
There exist constants $0<C_1^n\le C_2^n<\infty$ such that 
\begin{align}
\label{6.3}
C_1^n t^{-\frac{n}{2}+1}
\le \int_{\textbf{R}^n} e^{-t|\xi|^2}
\left|
\frac{\sin(t|\xi|)}{|\xi|}
\right|^2\,d\xi
\le C_2^n t^{-\frac{n}{2}+1}
\end{align}
for sufficiently large $t$. 
\end{lemma}

We omit the proof of Lemma~\ref{lem:6.3} 
in order to avoid repetition of the same calculation. 
Recall the proof of (i) in Lemma~\ref{lem:4.3} with $\gamma=0$, 
which is the original proof of Lemma~\ref{lem:6.3} given in \cite{I-2}. \\

\textbf{Proof of Theorem~\ref{thm:3.3}.} 
First, one can easily see that 
\begin{align*}
& \big\|
\hat{u}(t)
\big\|_2
\ge \big\|
\hat{u}(t)
\big\|_{L^2(|\xi|\le1)} 
\ge \biggr\|
E_1(t)\widehat{u_1}
\biggr\|_{L^2(|\xi|\le1)}
-\biggr\|
E_0(t)\widehat{u_0}
\biggr\|_{L^2(|\xi|\le1)}
-\biggr\|
\frac{|\xi|^2}{2}E_1(t)\widehat{u_0}
\biggr\|_{L^2(|\xi|\le1)} \\
& \ge \biggr\|
E_1(t)\widehat{u_1}
\biggr\|_{L^2(|\xi|\le1)}
-C\|u_0\|_1 (1+t)^{-\frac{n}{4}}
-C\|u_0\|_1 (1+t)^{-\frac{n}{4}-\frac{1}{2}} \\
& \ge \biggr\|
e_1^0(t) m[u_1]^0
\biggr\|_{L^2(|\xi|\le1)}
-\biggr\|
E_1(t)\widehat{u_1}
-e_1^0(t) 
m[u_1]^0
\biggr\|_{L^2(|\xi|\le1)} \\
& \qquad\qquad\qquad\qquad\quad
-C\|u_0\|_1 (1+t)^{-\frac{n}{4}}
-C\|u_0\|_1 (1+t)^{-\frac{n}{4}-\frac{1}{2}}.
\end{align*}
We have just used estimates [II] in Remark~\ref{rem:4.2} and \eqref{5.1}. 
Recalling \eqref{3.3} with $\gamma=0$, 
one has 
\[
\biggr\|
E_1(t)\widehat{u_1}
-e_1^0(t) 
m[u_1]^0
\biggr\|_{L^2(|\xi|\le1)}
=o(t^{-\frac{n}{4}+\frac{1}{2}}), 
\qquad
t\to\infty.
\]
Thus it follows that 
\begin{align*}
\big\|
\hat{u}(t)
\big\|_2
\ge \biggr\|
e_1^0(t) m[u_1]^0
\biggr\|_{L^2(|\xi|\le1)}
-o(t^{-\frac{n}{4}+\frac{1}{2}}),
\qquad
t\to\infty. 
\end{align*}
Rewriting
\[
e_1^0(t,\xi)
=e^{-\frac{t|\xi|^2}{2}}\frac{\sin(t|\xi|)}{|\xi|},
\qquad 
m[u_1]^0(\xi)
=\int_{\textbf{R}^n}u_1(x)\,dx, 
\]
the theorem follows from Lemmas~\ref{lem:6.1}, \ref{lem:6.2} and \ref{lem:6.3} since the same inequalities as \eqref{6.1}, \eqref{6.2} and \eqref{6.3} are also valid 
even if the integral domain is replaced 
from $\textbf{R}^n$ to $\{\xi\in\textbf{R}^n:|\xi|\le1\}$.
$\Box$
\\

\noindent
\textbf{Acknowledgment.}
The author would like to thank Professor Ryo Ikehata for useful discussions and warm encouragement.

\bibliographystyle{amsplain}

\end{document}